# ON THE ABSOLUTE CONTINUITY OF LÉVY PROCESSES WITH DRIFT

By Ivan Nourdin and Thomas Simon

*Université Henri Poincaré and Université d'Évry-Val d'Essonne*

We consider the problem of absolute continuity for the one-dimensional SDE

$$X_t = x + \int_0^t a(X_s)\,ds + Z_t,$$

where $Z$ is a real Lévy process without Brownian part and $a$ a function of class $\mathcal{C}^1$ with bounded derivative. Using an elementary stratification method, we show that if the drift $a$ is monotonous at the initial point $x$, then $X_t$ is absolutely continuous for every $t > 0$ if and only if $Z$ jumps infinitely often. This means that the drift term has a regularizing effect, since $Z_t$ itself may not have a density. We also prove that when $Z_t$ is absolutely continuous, then the same holds for $X_t$, in full generality on $a$ and at every fixed time $t$. These results are then extended to a larger class of elliptic jump processes, yielding an optimal criterion on the driving Poisson measure for their absolute continuity.

**1. Introduction and statement of the results.** Let $Z$ be a real Lévy process without Brownian part and jumping measure $\nu$, starting from 0. Without loss of generality, we may suppose that the sample paths of $Z$ are càdlàg. Consider the SDE

$$X_t = x + \int_0^t a(X_s)\,ds + Z_t, \tag{1}$$

where $a:\mathbb{R}\to\mathbb{R}$ is $\mathcal{C}^1$ with bounded derivative. Since $a$ is global Lipschitz, it is well known that there is a unique strong solution to (1). Let $\lambda$ be the Lebesgue measure on $\mathbb{R}$. For any real random variable $X$, we will write $X \ll \lambda$ for $X$ is absolutely continuous with respect to $\lambda$. We will say that a real function $f$ is increasing (resp. decreasing) at $x \in \mathbb{R}$ if there exists $\varepsilon > 0$ such that $f(y) < f(z)$ [resp. $f(y) > f(z)$] for every $x-\varepsilon < y < z < x+\varepsilon$, and that $f$ is monotonous at $x$ if it is either increasing or decreasing at $x$. The main purpose of this paper is to prove the following:

THEOREM A. *Suppose that $a$ is monotonous at $x$. Then the following equivalences hold:*

$$X_t \ll \lambda \quad \text{for every } t > 0 \iff X_1 \ll \lambda \iff \nu \text{ is infinite.}$$

This result is somewhat surprising because of the very weak optimal condition on the jumping measure $\nu$. Indeed, it has been known since Döblin (see Theorem 27.4 in [19]), that the following equivalence holds:

$$Z_t \text{ does not weight points for every } t > 0 \iff \nu \text{ is infinite.}$$

On the other hand, it is possible that $\nu$ is infinite and $Z_t$ is not absolutely continuous for every $t > 0$, see Theorem 27.19 in [19]. In other words, our result shows that a large class of drift perturbations may have some *regularizing effects* on the distribution of the perturbed Lévy process. It would be interesting to know if these effects also concern the regularity of the density.

Recall that the problem of absolute continuity for $Z_t$ itself is an irritating question, for which no necessary and sufficient condition on the jumping measure has been found as yet. As it had been pointed out by Orey (see Exercise 29.12 in [19]), a condition such that

$$\int_{|z|\le 1} |z|^\alpha \nu(dz) = +\infty$$

for some $\alpha \in ]0,2[$ is not sufficient because $\operatorname{Supp}\nu$ may become too sparse around zero. Up to now, the best criterion seems to be the following:

THEOREM (Kallenberg [11], Sato [18]). *Suppose that $\nu$ is infinite and let $\mu$ be the finite measure defined by $\mu(dx) = |x|^2(1+|x|^2)^{-1}\nu(dx)$. Then $Z_t$ is absolutely continuous for every $t > 0$ in the following situations:*

(a) *For some $n \ge 1$, the $n$th convolution power of $\mu$ is absolutely continuous.*
(b) $\lim_{\varepsilon\to 0} \varepsilon^{-2}|\log\varepsilon|^{-1}\mu(-\varepsilon,\varepsilon) = +\infty$.

The above condition (b) and the fact that $\nu$ is infinite when $\limsup_{\varepsilon\to 0}\varepsilon^{-2}\times \mu(-\varepsilon,\varepsilon) = +\infty$ and finite when $\lim_{\varepsilon\to 0}\varepsilon^{-2}|\log\varepsilon|^{r}\mu(-\varepsilon,\varepsilon) = 0$, for some $r > 1$, entail that the class of continuous, nonabsolutely continuous infinitely divisible distributions is actually very thin. However, Rubin and Tucker (see Theorem 27.23 in [19]) had constructed an exotic, atomic Lévy measure such that, for a given $t_0 > 0$, $Z_t$ is continuous singular for $t < t_0$ and absolutely continuous for $t > t_0$. This construction entails that absolute continuity is a *time-dependent* distributional property for Lévy processes, and we refer to the recent survey [22] for much more on this delicate topic, with connections to the notion of semi-self-decomposability and different families of algebraic numbers. Our result proves that there is no such transition phase phenomenon for a large class of drifted Lévy processes, in full generality on $\nu$.

We stress that the condition on $a$ is quite natural to obtain the equivalence stated in Theorem A. Suppose, indeed, that $a$ is constant and, say, positive, in an open neighborhood of $x$. Suppose that $Z$ has no drift, finite variation and infinitely many negative small jumps, but that $\operatorname{Supp}\nu$ is sparse enough in the neighborhood of zero for $Z_1$ to be only continuous singular. It is an easy consequence of a support theorem for jump processes in the finite variation case (see Theorem I in [21]) that, for every $\varepsilon > 0$, the event

$$\Omega_\varepsilon = \left\{ \sup_{0 \le t \le 1} |X_t - x| < \varepsilon \right\}$$

will have positive probability. Let now $A$ be a Borel set of Lebesgue measure zero such that $\mathbb{P}[Z_1 \in A] = 1$ and set $\tilde{A} = x + a(x) + A$. If $\varepsilon$ is chosen small enough, we have

$$\mathbb{P}[X_1 \in \tilde{A}] \ge \mathbb{P}[X_1 \in \tilde{A}, \Omega_\varepsilon] = \mathbb{P}[x + a(x) + Z_1 \in \tilde{A}, \Omega_\varepsilon] = \mathbb{P}[\Omega_\varepsilon] > 0,$$

which entails that $X_1$ is not absolutely continuous. For the same reason, Theorem A is not true in full generality when $a$ is locally flat on the left or on the right at the starting point $x$—see, however, Remark (a) below.

Equations of type (1) driven by Lévy processes are important for applications since they include storage processes, generalized OU processes, etc. They are also relevant in physics, for example, in climate models—see [9] and the references therein. Having in mind a good criterion of density for jumping SDEs, it is important to deal with a drift and a driving process which should be as general as possible. This is done in the next theorem, which yields the optimal level of generality for $Z$, and requires no further assumption on the drift $a$:

THEOREM B. *For every $t > 0$, the following inclusion holds:*

$$Z_t \ll \lambda \quad \Longrightarrow \quad X_t \ll \lambda.$$

Notice that this criterion for the absolute continuity of $X_t$ is the best possible on the driving process $Z$: when $a$ is constant, the inclusion is obviously an equivalence, and Theorem A depicts a situation where the reverse inclusion may not hold.

Over the last twenty years, a scattered literature dealing with the absolute continuity of SDEs with jumps has emerged, using various techniques—see, for example, [2, 4, 7, 8, 10, 17] and the older references therein. In the present paper we propose yet another way to tackle this question, which is mainly influenced by the "stratification method" due to Davydov, Lifshits and Smorodina [3]. The latter was introduced to study the absolute continuity of integral functionals of Brownian motion or more general Lévy processes, and relies roughly on a suitable decomposition of the underlying probability space into finite-dimensional strata. This method appears to be particularly well-adapted to our simple equation (1), and takes here an even more elementary form: we just need to choose *one* good jump of $Z$ with respect to which the drift term has a nonvanishing derivative, and then we use standard independence and distributional properties of Poisson measures.

In the next section we prove Theorems A and B. In the third and last section of this paper we generalize these results, without much effort, to the SDE

$$X_t = x + \int_0^t a(X_s)\, ds + \int_0^t \sigma(X_s) \diamond dZ_s,$$

where $\sigma$ is a $\mathcal{C}^1$ function with bounded derivative and $\diamond$ is the so-called Marcus integrator (see below for details) under the assumption that $\sigma$ never vanishes. This improves on a recent result of Ishikawa and Kunita [10], in the one-dimensional case.

## 2. Proofs.

2.1. *Proof of Theorem* A. As we said before, we first need to introduce a suitable decomposition of our underlying probability space. Let $\mathbf{F}$ be a closed set not containing 0 and such that $\nu(\mathbf{F}) \neq 0$. Set $\{T_n, n \geq 1\}$ for the sequence of jumping times of $\Delta Z$ into $\mathbf{F}$: $T_0 = 0$ and $T_n = \inf\{t > T_{n-1}/\Delta Z_t \in \mathbf{F}\}$ for every $n \geq 1$. Let $\mathcal{F}$ be the $\sigma$-algebra generated by $\{T_n, n \geq 2\}$, $\{\Delta Z_{T_n}, n \geq 1\}$ and the process $\tilde{Z}$ defined by

$$\tilde{Z}_t = Z_t - \sum_{T_n \leq t} \Delta Z_{T_n}$$

for every $t \geq 0$. Since $\{T_n, n \geq 1\}$ is the sequence of jumping times of some Poisson process independent of $\{\Delta Z_{T_n}, n \geq 1\}$ and $\tilde{Z}$, it is well known that $T_1$ has a uniform law on $[0, T_2]$ conditionally on $\mathcal{F}$. Hence, setting $T = T_1$ for simplicity, we can construct $Z$ on the disintegrated probability space $(\bar{\Omega} \times [0, T_2(\bar{\omega})], \mathcal{F} \times \mathcal{B}_{[0,T_2(\bar{\omega})]}, \bar{\mathbb{P}} \times \lambda^{\bar{\omega}}_{[0,T_2(\bar{\omega})]})$, where $\bar{Z}$ is the process defined by

$$\bar{Z}_t = Z_t - \mathbf{1}_{\{T \leq t\}} \Delta Z_T,$$

$(\bar{\Omega}, \mathcal{F}, \bar{\mathbb{P}})$ is the canonical space associated with $(\Delta Z_T, \bar{Z})$, and $\lambda^{\bar{\omega}}_{[0,T_2(\bar{\omega})]}$ is the normalized Lebesgue measure on $[0, T_2(\bar{\omega})]$. In other words, for every real functional $F$ of $Z$ and every Borelian $A$, we have the disintegration formula

$$\mathbb{P}[F \in A] = \mathbb{P}[F(\bar{\omega}, T) \in A] = \bar{\mathbb{E}}\left[\frac{1}{T_2(\bar{\omega})} \int_0^{T_2(\bar{\omega})} \mathbf{1}_{\{F(\bar{\omega},t) \in A\}}\, dt\right]. \tag{2}$$

We now come back to our equation (1). An important feature is that it can be solved *pathwise*, in the following sense: $X_t = Y_t + Z_t$ for every $t \geq 0$, where $\{Y_t, t \geq 0\}$ is the solution to the random ODE

$$Y_t = x + \int_0^t a(Y_s + Z_s)\, ds. \tag{3}$$

In the following we will work on $Y$ rather than on $X$. We begin with a fairly obvious lemma.

LEMMA 1. *Let $\rho : \mathbb{R}^+ \to \mathbb{R}$ be a càdlàg function and $a : \mathbb{R} \to \mathbb{R}$ be a $\mathcal{C}^1$ function with bounded derivative. Then the ODE*

$$x_t(x) = x + \int_0^t a(x_s(x) + \rho_s)\, ds$$

*has a unique solution which induces a differentiable flow $x \mapsto x_t(x)$. The derivative of the flow is given by*

$$\dot{x_t}(x) = 1 + \int_0^t \dot{a}(x_s(x) + \rho_s)\dot{x_s}(x)\, ds = \exp\left[\int_0^t \dot{a}(x_s(x) + \rho_s)\, ds\right].$$

PROOF. The existence of a unique solution to the ODE is a routine which follows from Picard's iteration scheme and Gronwall's lemma. The exponential formula for the derivative of the flow is not completely straightforward but, as for the case $\rho \equiv 0$, this can be obtained by introducing the ODE

$$u_t(x) = 1 + \int_0^t \dot{a}(x_s(x) + \rho_s) u_s(x)\, ds,$$

and proving that, for every $t > 0$ and $h \in \mathbb{R}$,

$$x_t(x+h) - x_t(x) - h u_t(x) = o(h)$$

as $h \to 0$. The latter estimation can be made in approximating the càdlàg sample paths of $\rho$ by continuous functions. We leave the details to the reader. □

We stress that, by a localization argument, we may (and will) suppose that $a$ is also bounded in order to prove the absolute continuity of $X_t$. We will note $M = \|a\|_\infty$. The following proposition yields a crucial computation, reminiscent of the Malliavin calculus.

PROPOSITION 2. *The map $T \mapsto Y_1(\bar{\omega}, T)$ is everywhere differentiable except possibly on $T = 1$, with derivative given by*

$$\frac{dY_1}{dT} = \mathbf{1}_{\{T<1\}}(a(X_{T-}) - a(X_T)) \exp\left[\int_T^1 \dot{a}(X_s)\,ds\right]$$

*when $T \neq 1$.*

PROOF. It is obvious that $Y_1$ does not depend on $T$ as soon as $T > 1$. Fix now $T \in ]0,1[$. For any $h > 0$ small enough, define the process

$$Z_t^h = Z_t + \Delta Z_T(\mathbf{1}_{\{T-h\le t\}} - \mathbf{1}_{\{T\le t\}})$$

and set $Y^h$ for the solution to (3), with parameter $Z^h$ instead of $Z$. Notice first that

$$Y_1 - Y_1^h = \left(Y_T + \int_T^1 a(Y_s + Z_s)\,ds\right) - \left(Y_T^h + \int_T^1 a(Y_s^h + Z_s)\,ds\right). \tag{4}$$

Besides,

$$\begin{aligned} Y_T - Y_T^h = & \int_{T-h}^T (a(Y_s + Z_s) - a(Y_s + Z_s^h))\,ds \\ & + \int_{T-h}^T (a(Y_s + Z_s^h) - a(Y_s^h + Z_s^h))\,ds. \end{aligned}$$

By the boundedness of $a$ and the fact that $Y_{T-h} = Y_{T-h}^h$, we see that there exists $M_1 > 0$ such that $|Y_s - Y_s^h| \le M_1 h$ for every $s \in [T-h, T]$. Since $\dot{a}$ is also bounded, this entails that

$$Y_T - Y_T^h = \int_{T-h}^T (a(Y_s + Z_s) - a(Y_s + Z_s^h))\,ds + O(h^2).$$

Hence,

$$\begin{aligned} \frac{Y_T - Y_T^h}{h} & = \left(\frac{1}{h}\int_{T-h}^T a(X_s)\,ds\right) - \left(\frac{1}{h}\int_{T-h}^T a(X_s + \Delta X_T)\,ds\right) + O(h) \\ & \longrightarrow a(X_{T-}) - a(X_T) \qquad \text{as } h \downarrow 0, \end{aligned} \tag{5}$$

where the convergence comes from the fact that $s \to X_s$ has left-hand limits. It follows from (4), (5) and Lemma 1 applied to the ODE

$$Y_{T+t} = x_t(Y_T) = Y_T + \int_0^t a(x_s(Y_T) + Z_{T+s})\,ds$$

that

$$\lim_{h\downarrow 0}\left(\frac{Y_1 - Y_1^h}{h}\right) = (a(X_{T-}) - a(X_T)) \exp\left[\int_T^1 \dot{a}(X_s)\,ds\right]. \tag{6}$$

Similarly, defining

$$Z_t^h = Z_t + \Delta Z_T(\mathbf{1}_{\{T+h\leq t\}} - \mathbf{1}_{\{T\leq t\}})$$

for $h$ small enough, and setting $Y^h$ for the solution to (3) with parameter $Z^h$ instead of $Z$, we can prove that

$$\lim_{h\downarrow 0}\left(\frac{Y_1^h - Y_1}{h}\right) = (a(X_{T-}) - a(X_T))\exp\left[\int_T^1 \dot{a}(X_s)\,ds\right]. \tag{7}$$

Putting (6) and (7) together completes the proof in the case $T\neq 1$. Last, we notice that $T\mapsto Y_1(\bar{\omega},T)$ has obviously a right derivative at $T=1$ which is zero and, reasoning exactly as above, a left derivative at $T=1$ which is $(a(X_{T-}) - a(X_T))$ and may be non zero. □

Proof of Theorem A. The first inclusion is trivial, and the second is easy: if $\nu$ is finite, then $Z$ does not jump up to time 1 with positive probability, so that if $d$ stands for the drift of $Z$, then $X_1$ has an atom at $x_1$ where $\{x_t, t\geq 0\}$ solves the equation

$$x_t = x + \int_0^t a(x_s)\,ds + \mathrm{d}t.$$

It remains to show that if $\nu$ is infinite, then $X_t \ll \lambda$ for every $t>0$. Considering $-X$ if necessary, we may suppose that $\nu([0,1]) = +\infty$ and that $a$ is increasing on $]x-\varepsilon, x+\varepsilon[$ for some $\varepsilon>0$, and we may fix $t=1$. Let $A$ be a Borelian set of $\mathbb{R}$ such that $\lambda(A)=0$. We need to prove that, for every $\delta>0$,

$$\mathbb{P}[X_1\in A] < \delta. \tag{8}$$

Fix $\delta>0$. By the right-continuity of $Z$, there exists $\beta>0$ such that

$$\mathbb{P}\left[\sup_{s\leq\beta}|Z_s|\geq\varepsilon/6\right] < \delta/2.$$

Set $\gamma = \beta\wedge(\varepsilon/3M)\wedge 1$ (recall that $M=\|a\|_\infty$) and $\{T_\varepsilon^\eta(n), n\geq 1\}$ for the sequence of jumping times of $\Delta Z$ into $\{\eta\leq z\leq\varepsilon/6\}$, with $\eta$ chosen small enough. Since $\nu([0,1]) = +\infty$, there exists $\eta$ such that

$$\mathbb{P}[T_\varepsilon^\eta(2)\geq\gamma] < \delta/2.$$

We note $T = T_\varepsilon^\eta(1)$, $T_2 = T_\varepsilon^\eta(2)$, and $\bar{Z}_t = Z_t - \Delta Z_T\mathbf{1}_{\{T\leq t\}}$ for every $t\geq 0$. Let $\mathcal{F}$ be the $\sigma$-algebra generated by $\Delta Z_T$ and the process $\bar{Z}$. Conditionally on $\mathcal{F}$, $T$ has uniform law on $[0,T_2]$ and with the same notation as above, we can work on the disintegrated probability space $(\bar{\Omega}\times[0,T_2(\bar{\omega})], \mathcal{F}\times\mathcal{B}_{[0,T_2(\bar{\omega})]}, \bar{\mathbb{P}}\times\lambda_{[0,T_2(\bar{\omega})]})$. Notice that if we set

$$\bar{\Omega}_1 \stackrel{\text{def}}{=} \left\{\sup_{s\leq\gamma}|\bar{Z}_s| < \varepsilon/3, T_\varepsilon^\eta(2) < \gamma\right\} \supset \left\{\sup_{s\leq\gamma}|Z_s| < \varepsilon/6, T_\varepsilon^\eta(2)<\gamma\right\},$$

then $\bar{\Omega}_1 \in \mathcal{F}$ and $\mathbb{P}[\bar{\Omega}_1^c] < \delta$. Hence,

$$\mathbb{P}[X_1 \in A] < \delta + \mathbb{E}[\mathbf{1}_{\bar{\Omega}_1} \mathbb{P}[X_1 \in A|\mathcal{F}]].$$

We will now prove that

$$\mathbb{P}[X_1 \in A|\mathcal{F}](\bar{\omega}) = 0,$$

for every $\bar{\omega} \in \bar{\Omega}_1$, which will yield (8) and complete the proof of the theorem. Fix $\bar{\omega} \in \bar{\Omega}_1$ once and for all. The key-point is that $Z_1(\bar{\omega}, T)$ does *not* depend anymore on $T$: indeed, we know that the process $Z(\bar{\omega}, T)$ jumps at least once into $\{\eta \le z \le \varepsilon/6\}$ before time 1, and by the Lévy–Itô decomposition, the terminal value $Z_1$ is independent of the first jumping time $T$. Hence, we can write

$$\begin{aligned}\mathbb{P}[X_1 \in A|\mathcal{F}](\bar{\omega}) &= \mathbb{P}[Y_1(\bar{\omega}, T) \in A - Z_1(\bar{\omega})|\mathcal{F}](\bar{\omega}) \\ &= \frac{1}{T_2(\bar{\omega})} \int_0^{T_2(\bar{\omega})} \mathbf{1}_{\{Y_1(\bar{\omega}, t) \in A - Z_1(\bar{\omega})\}}\, dt,\end{aligned}$$

where the last equality comes from (2). Besides, again since $T < T_2(\bar{\omega}) < 1$, it follows from Proposition 2 that

$$\frac{dY_1}{dT}(\bar{\omega}, T) = (a(X_{T-}(\bar{\omega})) - a(X_T(\bar{\omega}))) \exp\left[\int_T^1 \dot{a}(X_s(\bar{\omega}, T))\, ds\right].$$

Another important point is that, for every $T \in ]0, T_2(\bar{\omega})[$,

$$|X_T(\bar{\omega}) - x| \le |\bar{Z}_T(\bar{\omega})| + |\Delta Z_T(\bar{\omega})| + TM < \varepsilon/3 + \varepsilon/6 + \varepsilon/3 < \varepsilon$$

and

$$|X_{T-}(\bar{\omega}) - x| \le |\bar{Z}_T(\bar{\omega})| + TM < \varepsilon/3 + \varepsilon/3 < \varepsilon.$$

Hence, since $a$ is increasing on $]x - \varepsilon, x + \varepsilon[$ and $X_T(\bar{\omega}) - X_{T-}(\bar{\omega}) = \Delta Z_T \ge \eta > 0$, we get

$$\frac{dY_1}{dT}(\bar{\omega}, T) < 0$$

for every $T \in ]0, T_2(\bar{\omega})[$, so that $T \mapsto Y_1(\bar{\omega}, T)$ is a diffeomorphism on $]0, T_2(\bar{\omega})[$. But $A - Z_1(\bar{\omega})$ has Lebesgue measure zero, and this yields finally

$$\int_0^{T_2(\bar{\omega})} \mathbf{1}_{\{Y_1(\bar{\omega}, t) \in A - Z_1(\bar{\omega})\}}\, dt = 0,$$

as desired. $\square$

REMARKS. (a) It is easy to see that the proof of Theorem A carries over to two other particular situations: when $a(x) > 0$ [resp. $a(x) < 0$], $a$ is monotonous on the right (resp. on the left) of $x$ and $Z$ (resp. $-Z$) is a subordinator. The case when $x$ is an isolated point of $\{y/\dot{a}(y)=0\}$ could also probably be handled in the same way, because a.s. $X_t$ visits instantaneously either $]-\infty, x[$ or $]x, +\infty[$. However, thinking of the case where $\dot{a}$ varies very irregularly around $x$, it seems difficult to relax further the local monotonicity condition on $a$ without asking more from the driving process $Z$.

(b) When $Z_1$ is continuous singular, Theorem A entails that the laws of the processes $\{Z_t, t \geq 0\}$ and $\{X_t, t \geq 0\}$ are mutually singular under the monotonicity assumption on $a$. When $a(x) \neq 0$ and $Z$ has finite variation, this is also a simple consequence of the a.s. fact that

$$\lim_{t \downarrow 0} t^{-1} X_t = a(x) + \lim_{t \downarrow 0} t^{-1} Z_t,$$

which is a consequence of Theorem 43.20 in [19] and a comparison argument. When $a \equiv k \neq 0$, this mutual singularity is also a consequence of Skorohod's dichotomy theorem—see Theorem 33.1 in [19], or [13] for a recent extension of the latter result to general transformations on Poisson measures. We got stuck in proving this mutual singularity under the sole assumptions that $x + Z_t$ visits a.s. $\{y/a(y) \neq 0\}$ and that $Z$ has no Gaussian part.

2.2. *Proof of Theorem* B. Without loss of generality, we can suppose $t = 1$. Let $\{\mathcal{F}_t, t \geq 0\}$ be the completed natural filtration of $Z$. For every $c > 0$, we define the closed set $\mathbf{D}_c = \{x \in \mathbb{R}/|\dot{a}(x)| \geq c\}$ and the $\mathcal{F}_t$-stopping time $T_c = \inf\{t \geq 0/X_t \in \mathbf{D}_c\}$. Let $A$ have zero Lebesgue measure. Since $T_c$ is measurable with respect to $\mathcal{F}_{T_c}$, we can decompose

$$\mathbb{P}[X_1 \in A] = \mathbb{P}[T_c \geq 1, X_1 \in A] + \mathbb{P}[T_c < 1, \mathbb{E}[X_1 \in A/\mathcal{F}_{T_c}]].$$

By the strong Markov property for $X$, we have a.s. on $\{T_c < 1\}$

$$\mathbb{E}[X_1 \in A/\mathcal{F}_{T_c}] = \mathbb{E}[\tilde{X}_{1-T_c} \in A/\mathcal{F}_{T_c}],$$

where $\tilde{X}$ solves the SDE

$$\tilde{X}_t = X_{T_c} + \int_0^t a(\tilde{X}_s)\, ds + \tilde{Z}_t,$$

with $\tilde{Z}$ a copy of $Z$ independent of $\mathcal{F}_{T_c}$. Recall that $(T_c, X_{T_c})$ is measurable with respect to $\mathcal{F}_{T_c}$ and that $|\dot{a}(X_{T_c})| \geq c > 0$ by right-continuity. Hence, since the jumping measure of $\tilde{Z}$ is infinite, we can apply Theorem A and get

$$\mathbb{E}[\tilde{X}_{1-T_c} \in A/\mathcal{F}_{T_c}] = 0,$$

a.s. on $\{T_c < 1\}$. In particular, $\mathbb{P}[X_1 \in A] = \mathbb{P}[T_c \geq 1, X_1 \in A]$ for every $c > 0$. By right-continuity, $T_c \downarrow T_0 = \inf\{t \geq 0/\dot{a}(X_t) \neq 0\}$ as $c \downarrow 0$, so that finally

$$\mathbb{P}[X_1 \in A] = \mathbb{P}[T_0 \geq 1, X_1 \in A].$$

On the event $\{T_0 \geq 1\}$, it follows readily from Itô's formula with jumps that, for every $t \in [0,1]$,

$$a(X_t) = a(x) + \sum_{s \leq t} (a(X_s) - a(X_{s-})). \tag{9}$$

For every $\eta > 0$, introduce the events $\Omega_a^{\eta+} = \{\exists t \in ]0,1[/a(X_t) - a(X_{t-}) \geq \eta\}$, $\Omega_a^{\eta-} = \{\exists t \in ]0,1[/a(X_t) - a(X_{t-}) \leq -\eta\}$, and set

$$\Omega_a = \bigcup_{\eta > 0} (\Omega_a^{\eta+} \cup \Omega_a^{\eta-}) = \lim_{\eta \downarrow 0} \uparrow (\Omega_a^{\eta+} \cup \Omega_a^{\eta-}).$$

On $\Omega_a^c \cap \{T_0 \geq 1\}$, it follows readily from (9) that a.s. $X_1 = a(x) + Z_1$, whence

$$\mathbb{P}[\Omega_a^c \cap \{T_0 \geq 1\}, X_1 \in A] \leq \mathbb{P}[Z_1 \in A - a(x)] = 0,$$

from the assumption on $Z$. It remains to consider the last situation where $\dot{a}(X) \equiv 0$ but the process $a(X)$ jumps in the interval $[0,1]$, and we are now going to show that, for every $\delta, \eta > 0$,

$$\mathbb{P}[\Omega_a^{\eta-}, T_0 \geq 1, X_1 \in A] + \mathbb{P}[\Omega_a^{\eta+}, T_0 \geq 1, X_1 \in A] < \delta, \tag{10}$$

which will give $\mathbb{P}[X_1 \in A] < \delta$ from what precedes and, hence, conclude the proof of Theorem B. The argument is very similar to that of Theorem A and we will only use the infiniteness of $\nu$. We first consider the event $\Omega_a^{\eta+} \cap \{T_0 \geq 1\} \cap \{X_1 \in A\}$. Fix $\delta, \eta > 0$. Let $S_\eta$ denote the first jump of $a(X)$ greater than $\eta$ and, for every $\varepsilon > 0$, let $\{T_n^\varepsilon, n \geq 1\}$ be the ordered sequence of jumping times of $Z$ into $\{|z| \geq \varepsilon\}$. Since $a$ is global Lipschitz and $\Delta X_t = \Delta Z_t$ for every $t \geq 0$, we can take $\varepsilon > 0$ small enough such that $S_\eta \in \{T_n^\varepsilon, n \geq 1\}$ a.s. For every $\varepsilon > 0$, consider now the event

$$\Omega_\varepsilon = \{\exists n \geq 2/S_\eta < T_n^\varepsilon < 1\}.$$

Since $\nu$ is infinite, by density we can choose $\varepsilon$ small enough such that

$$\mathbb{P}[\Omega_a^{\eta+} \cap \Omega_\varepsilon^c] < \delta/2. \tag{11}$$

On the other hand,

$$\begin{aligned} &\mathbb{P}[\Omega_a^{\eta+} \cap \Omega_\varepsilon, T_0 \geq 1, X_1 \in A] \\ &\quad = \sum_{n=1}^{+\infty} \mathbb{P}[\Omega_a^{\eta+}, S_\eta = T_n^\varepsilon, T_{n+1}^\varepsilon < 1 \leq T_0, X_1 \in A]. \end{aligned} \tag{12}$$

For every $n \geq 1$, let $\mathcal{F}_n$ be the $\sigma$-algebra generated by the process

$$Z_t^n = Z_t - \mathbf{1}_{\{T_n^\varepsilon \leq t\}} \Delta Z_{T_n^\varepsilon}$$

and the random variable $\Delta Z_{T_n^\varepsilon}$. Notice the a.s. inclusion

$$\begin{aligned} &\{\Omega_a^{\eta+}, S_\eta = T_n^\varepsilon, T_{n+1}^\varepsilon < 1 \leq T_0, X_1 \in A\} \\ &\qquad \subset \{a(X_{T_n^\varepsilon}) - a(X_{T_n^\varepsilon -}) \geq \eta, T_{n+1}^\varepsilon < 1 \leq T_0, X_1 \in A\} \end{aligned}$$

and recall that, conditionally on $\mathcal{F}_n$, $T_n^\varepsilon$ has a uniform law on $[T_{n-1}^\varepsilon, T_{n+1}^\varepsilon]$. Hence, using the same notation and arguments as in Theorem A, we have

$$\begin{aligned}
&\mathbb{P}[\Omega_a^{\eta+}, S_\eta = T_n^\varepsilon, T_{n+1}^\varepsilon < 1 \le T_0, X_1 \in A | \mathcal{F}_n](\bar\omega) \\
&\quad \le \frac{\mathbf{1}_{\{T_{n+1}^\varepsilon(\bar\omega) < 1\}}}{(T_{n+1}^\varepsilon(\bar\omega) - T_{n-1}^\varepsilon(\bar\omega))} \\
&\qquad \times \int_{T_{n-1}^\varepsilon(\bar\omega)}^{T_{n+1}^\varepsilon(\bar\omega)} \mathbf{1}_{\{a(X_t(\bar\omega)) - a(X_{t-}(\bar\omega)) \ge \eta, T_0(\bar\omega,t) \ge 1, Y_1(\bar\omega,t) \in A - Z_1(\bar\omega)\}}\, dt = 0,
\end{aligned}$$

where in the second equality we used the fact, proved in Proposition 2, that, under the integral,

$$\begin{aligned}
\frac{dY_1}{dt}(\bar\omega, t) &= (a(X_{t-}(\bar\omega)) - a(X_t(\bar\omega))) \exp\left[\int_t^1 \dot a(X_s(\bar\omega,t))\, ds\right] \\
&= a(X_{t-}(\bar\omega)) - a(X_t(\bar\omega)) \\
&\le -\eta
\end{aligned}$$

for every $t \in [T_{n-1}^\varepsilon, T_{n+1}^\varepsilon]$. Hence, from (11) and (12), we get

$$\mathbb{P}[\Omega_a^{\eta+}, T_0 \ge 1, X_1 \in A] < \delta/2.$$

We can prove $\mathbb{P}[\Omega_a^{\eta-}, T_0 \ge 1, X_1 \in A] < \delta/2$ similarly, which yields (10) as desired.

Remark. A quicker proof of Theorem B could be the following: by a theorem of Sharpe (see Theorem 3.3 in [20]) we know that when $Z_t$ has a density, then the latter is positive all over the interior of its support which is either $\mathbb{R}$ itself, or a half-line. Besides, including the drift coefficient of $Z$ in the function $a$ if necessary, we can suppose that this half-line is $\mathbb{R}^+$ or $\mathbb{R}^-$. Hence, for almost every $z \in \operatorname{Supp} Z_1$, it is possible to define $\mathbb{P}_1^z$ the law of the Lévy bridge associated with $Z$ from 0 to $z$ in time 1. Let $\mu$ (resp. $\nu$) now be the absolutely continuous (resp. singular) part of the law of $X_1$, and let $A$ be a Borel set with null Lebesgue measure such that $\nu(A) = \nu(\mathbb{R})$. We have

$$\mathbb{P}[X_1 \in A] = \int_{\mathbb{R}} \mathbb{P}_1^z[Y_1 \in A - z] p_{Z_1}(z)\, dz.$$

On the one hand, it follows easily from the independence of the increments of $Z$, the absolute continuity of $Z_2$ and the fact that $\operatorname{Supp} Z_1 = \operatorname{Supp} Z_2$, that $\mathbb{P}[Y_1 \in A - z] = 0$ a.s. on $z \in \operatorname{Supp} Z_1$. On the other hand, we can suppose by standard approximation that $a$ has compact support, so that $t \mapsto Y_t$ is a global Lipschitz process. In particular, for some $\kappa > 0$, setting $A_{\kappa t} = \{x / d(x, A) < \kappa(1-t)\}$, we see that $\mathbb{P}[Y_t \in A_{\kappa t} - z]$ *decreases* to $\mathbb{P}[Y_1 \in A - z]$ as $t \uparrow 1$. By equivalence of $\mathbb{P}_1^z$ and $\mathbb{P}$ on $\sigma\{Z_s, s \le t\}$, we have $\mathbb{P}_1^z[Y_t \in A_{\kappa t} -$ $z] < 1$ and by monotonicity, $\mathbb{P}_1^z[Y_t \in A - z] < 1$. This holding a.s. on $z \in \operatorname{Supp} Z_1$, we have $\nu(\mathbb{R}) = \mathbb{P}[X_1 \in A] < 1$, so that, in particular, $\mu \not\equiv 0$. Now since $X_1$ is the solution to a regular SDE defining a good flow of diffeomorphisms [12], it is hard to believe that $\nu$ should not be zero. However, we could not prove this last fact. In particular, it seems hard to transpose in this Markovian but non Lévy framework the classical argument of Hartman and Wintner (see Theorem 27.16 in [19]) for pureness of discretely generated ID distributions.

**3. Extension to some elliptic jump processes.** We now wish to extend the above results to a more general class of SDEs driven by a one-dimensional Lévy process, the so-called Marcus equations [14]. More precisely, we consider

$$X_t = x + \int_0^t a(X_s)\,ds + \int_0^t \sigma(X_s) \diamond dZ_s, \tag{13}$$

where $a$ and $\sigma$ are $\mathcal{C}^1$ with bounded derivative. The stochastic integral is defined in the following way:

$$\int_0^t \sigma(X_s) \diamond dZ_s = \int_0^t \sigma(X_{s-}) \circ dZ_s + \sum_{s \le t} \rho(X_{s-}, \Delta Z_s),$$

where $\circ$ is the classical Meyer–Stratonovitch integrator and $\rho : \mathbb{R}^2 \to \mathbb{R}$ is a local Lipschitz function such that, when $(y, z)$ stays in a fixed compact set of $\mathbb{R}^2$,

$$|\rho(y,z)| \le K|z|^2 \tag{14}$$

for some constant $K$, and such that in (13), each time $t$ when $Z$ jumps, $X_t$ is given by the integral in time 1 of the vector field $x \mapsto \sigma(x)\Delta Z_t$, starting from $X_{t-}$. With this definition of $\diamond$, it is well known that (13) has a unique strong solution under the above conditions on $a$ and $\sigma$.

A thorough study of Marcus equations driven by general (multidimensional) càdlàg semimartingales can be found in [12]. Notice that, in the literature, Marcus equations are sometimes called "canonical" equations; see [10] and the references therein. These Marcus equations have the major drawback that they only concern a specific class of integrands. Indeed, classical stochastic calculus allows us to define SDEs of type (13) where $\rho$ is *any* function verifying (14). This larger class of jump functions yields the required level of generality for applications. However, Marcus' choice of $\rho$ makes it possible to get a first-order change of variable formula, which simplifies the computations quite considerably. More precisely, we have

$$\begin{aligned} f(X_t) &= f(x_0) + \int_0^t f'(X_s) \diamond dX_s \\ &= f(x_0) + \int_0^t f'a(X_s)\,ds + \int_0^t f'\sigma(X_s) \diamond dZ_s \end{aligned} \tag{15}$$

for smooth real functions $f$—see Section 4 in [12] for details. Beware that, in (15), the jumps of $f(X_t)$ are defined in a different manner than those of $X_t$ in (13), because of the integro-differential term in Itô's formula with jumps. However, it follows from Definition 4.1. in [12] that when $f'\sigma \equiv k$ a constant function, then

$$f(X_t) = f(x_0) + \int_0^t f'a(X_s)\,ds + kZ_t,$$

a fact which will be used subsequently. It would be interesting to know if Marcus' choice of $\rho$ is the *only* one yielding a formula such as (15).

We will make further on the following *uniform ellipticity assumption* on the diffusion coefficient $\sigma$:

ASSUMPTION H. The function $\sigma$ does not vanish.

In order to simplify the presentation, we first consider the case when $Z$ has no Brownian part. Introducing the function $b \equiv a/\sigma$, the following corollary is an easy consequence of Theorems A and B.

COROLLARY. *Assume that $Z$ has no Brownian part and that Assumption H holds. If $b$ is monotonous in an open neighborhood of $x$, then we have the following equivalences:*

$$X_t \ll \lambda \quad \text{for every } t > 0 \iff X_1 \ll \lambda \iff \nu \text{ is infinite.}$$

*Moreover, for every $t > 0$, we have the inclusion*

$$Z_t \ll \lambda \implies X_t \ll \lambda.$$

PROOF. Define

$$f : x \mapsto \int_0^x \frac{dt}{\sigma(t)}$$

and set $Y_t = f(X_t)$ for every $t \ge 0$. Since $\sigma$ is continuous and never vanishes, $f$ is a diffeomorphism and it is clear that

$$X_t \ll \lambda \iff Y_t \ll \lambda$$

for every $t > 0$. Itô's formula (15) yields

$$\begin{aligned} Y_t &= f(x) + \int_0^t f'(X_s) \diamond dX_s \\ &= f(x) + \int_0^t \left(\frac{a}{\sigma}\right)(X_s)\,ds + \int_0^t \left(\frac{\sigma}{\sigma}\right)(X_s) \diamond dZ_s \qquad (16) \\ &= f(x) + \int_0^t b \circ f^{-1}(Y_s)\,ds + Z_t. \end{aligned}$$

When $b$ is monotonous in an open neighborhood of $x$, then $b \circ f^{-1}$ is monotonous in an open neighborhood of $f(x)$, since $f$ is a diffeomorphism. Hence, from (16), we just need to apply Theorem A to obtain the equivalences. When $b$ is not monotonous, Theorem B yields the desired inclusion. □

REMARK. The above corollary shows, for example, that $X_1$ is absolutely continuous as soon as $\nu$ is infinite, Assumption H holds and $[a, \sigma](x) \neq 0$, where $[\cdot, \cdot]$ stands for the standard Lie bracket. The latter condition may be viewed as a restricted Hörmander condition. Actually, in the multidimensional framework, the absolute continuity of $X_1$ was recently investigated under various Hörmander conditions [10]—see also [17] for a larger class of Itô equations. However, these papers (as well as, to our knowledge, all articles on Malliavin's calculus with jumps) require additional assumptions on $Z$ which are much stronger than the sole infiniteness of $\nu$ or the absolute continuity of $Z_1$. On the other hand, our class of equations is a bit restrictive when compared to the references quoted in the Introduction. When $\nu$ itself is absolutely continuous, one can use the simple method of Fournier and Giet [8] to prove that $X_t$ has a density for every $t > 0$. Similarly, the more involved method of Picard [17] applies when $\nu$ is singular but has enough small jumps. Both references deal with more general equations than (13). Our main motivation to write this paper was to get an *optimal* criterion on the driving process, and we do not know as yet what this criterion should be in a more general (e.g., multidimensional) framework.

In the following proposition, which is fairly straightforward and independent of what precedes, we obtain another necessary and sufficient condition for the absolute continuity of $Z_t$, in the particular situation when $a$ and $\sigma$ are proportional.

PROPOSITION 3. *Assume Assumption* H *and that* $a \equiv k\sigma$ *for some* $k \in \mathbb{R}$. *Then, for every* $t > 0$,

$$X_t \ll \lambda \quad \Longleftrightarrow \quad Z_t \ll \lambda.$$

PROOF. We can write

$$X_t = x + \int_0^t \sigma(X_s) \diamond dZ_s^k,$$

where $Z^k$ stands for the linearly drifted process: $Z_t^k = Z_t + kt$ for every $t \geq 0$. Let $(y, t) \mapsto \varphi(y, t)$ be the flow associated to $\sigma$: $\dot{\varphi}_t(y, t) = \sigma(\varphi(y, t))$, $\varphi(y, 0) = y$ for every $(y, t) \in \mathbb{R}^2$. Since $\sigma$ does not vanish, we see that $t \mapsto$ $\varphi(y,t)$ is a diffeomorphism. Besides, it follows immediately from Itô's formula (15) and the unicity of solutions to (13) that $X_t = \varphi(x, Z_t^k)$ for every $t \geq 0$. We get finally

$$X_t \ll \lambda \quad \iff \quad Z_t^k \ll \lambda \quad \iff \quad Z_t \ll \lambda$$

for every $t > 0$. $\square$

REMARK. Thinking of the equation $dX_t = X_t \diamond dZ_t$, whose solution is $X_t = X_0 \exp Z_t$, it is clear that Assumption H is not necessary in the statement of the above proposition.

To conclude this paper, and for the sake of completeness, we consider the easy situation when $Z$ has a Brownian part. Notice that here no assumption on $\nu$ is required.

PROPOSITION 4. *Assume Assumption* H *and that* $Z$ *has a nontrivial Brownian part. Then* $X_t \ll \lambda$ *for every* $t > 0$.

PROOF. When $a \equiv 0$, then we can reason exactly as in Proposition 3, since, by convolution, $Z_t$ has obviously a $\mathcal{C}^\infty$ density. When $a \not\equiv 0$, we can first suppose by an approximation argument that $a$ has compact support. We then use the Doss–Sussman transformation, which was established in [5] for Marcus equations on $\mathbb{R}$: for every $t > 0$, we write

$$X_t = \varphi(Y_t, Z_t),$$

where $\varphi$ is defined as in the proof of Proposition 3 and $Y$ is the solution to the random ODE

$$Y_t = x_0 + \int_0^t b(Y_s, Z_s)\, ds,$$

with the notation

$$b(x,y) = \frac{a(\varphi(x,y))}{\dot{\varphi}_x(x,y)} = a(\varphi(x,y)) \exp - \left[ \int_0^y \dot{\sigma}(\varphi(a,u))\, du \right] \tag{17}$$

for every $(x,y) \in \mathbb{R}^2$. We first notice that, for every fixed $x \in \mathbb{R}$, Assumption H entails that $t \mapsto \varphi(x,t)$ is a bijection onto $\mathbb{R}$. Indeed, without loss of generality, we can suppose that $\sigma > 0$, so that $\varphi(x,\cdot)$ is increasing and $\ell = \lim_{t\to+\infty} \varphi(x,t)$ exists in $\mathbb{R} \cup \{+\infty\}$. If $\ell \neq +\infty$, then $\lim_{t\to\infty} \varphi_t'(x,t) = \sigma(\ell) > 0$ and $\lim_{t\to+\infty} \varphi(x,t) = +\infty$. Similarly, we can show that $\lim_{t\to-\infty} \varphi(x,t) = -\infty$, which yields the desired bijection property. We will denote by $\psi : \mathbb{R}^2 \to \mathbb{R}$ the inverse of $\varphi$, that is, for each fixed $x \in \mathbb{R}$, $\psi(x,t)$ is the unique solution to $\varphi(x, \psi(x,t)) = t$. Using the flow property of $\varphi$, we get

$$X_t = \varphi(\varphi(x, \psi(x, Y_t)), Z_t) = \varphi(x, Z_t + \psi(x, Y_t)),$$

so that, by injectivity, it suffices to prove that $Z_t + \psi(x, Y_t)$ itself is absolutely continuous. Since $a$ has compact support, it is easy to see from (17) and the bijection property of $\varphi$ that $Y_t$ is a bounded random variable for every $t > 0$, and clearly $\psi(x, Y_t)$ is also bounded for every $t > 0$. Fix now $A$ a set of Lebesgue measure 0. Write $\mathbb{P} = \mathbb{P}^1 \otimes \mathbb{P}^2$ where $\mathbb{P}^1$ stands for the Poissonian part of $Z$ and $\mathbb{P}^2$ for its nontrivial Brownian part, and decompose $Z$ into $Z(\omega) = Z^1(\omega_1) + Z^2(\omega_2)$ accordingly: $Z^1$ is a Lévy process without Brownian part and $Z^2$ a rescaled Brownian motion. Since $t \mapsto \psi(x, Y_t)(\omega_1, .)$ is bounded continuous for almost every $\omega_1$, we can apply Girsanov's theorem under the measure $\mathbb{P}^2$ and get

$$\begin{aligned}\mathbb{P}[Z_t + \psi(x, Y_t) \in A] &= \mathbb{P}^1 \otimes \mathbb{P}^2[Z_t^2(\omega_2) + \psi(x, Y_t)(\omega_1, \omega_2) \in A - Z_t^1(\omega_1)] \\ &= \int d\mathbb{P}^1(\omega_1)\tilde{\mathbb{P}}^2_{\omega_1}[Z_t^2(\omega_2) \in A - Z_t^1(\omega_1)],\end{aligned}$$

where $\tilde{\mathbb{P}}^2_{\omega_1}$ is equivalent to $\mathbb{P}^2$ for almost every $\omega_1$. Since $A - Z_t^1(\omega_1)$ has zero Lebesgue measure, we get $\tilde{\mathbb{P}}^2_{\omega_1}[Z_t^2(\omega_2) \in A - Z_t^1(\omega_1)] = 0$ for almost every $\omega_1$, so that

$$\mathbb{P}[Z_t + \psi(x, Y_t) \in A] = 0,$$

as desired. □

Final remarks. (a) It is well known that Assumption H is far from being necessary when $Z$ is Brownian motion, according to the classical Bouleau–Hirsch criterion [1] which says that, in this case, $X_t$ has a density if and only if $t > t_0$, where $t_0$ is the first (deterministic) entrance-time into $\{\sigma(X_s) \neq 0\}$. It is somewhat tantalizing, but maybe challenging, to try to generalize this criterion to Marcus equations driven by a general one-dimensional Lévy process. In this direction, let us mention that Coquio [2] had proved a "density-image property for the energy measure" in a Poissonian framework.

(b) In the companion paper [16], we prove an analogue of Theorem A when the driving process is a fractional Brownian motion of any Hurst index, and we also extend Bouleau–Hirsch's criterion in this framework.

(c) At the beginning of this research, we had first tried to prove Theorem B and the corresponding corollary in the following way: assuming, without loss of generality, that $a$ is bounded, let $k > 0$ be such that $a(x) < k - 1$ for every $x \in \mathbb{R}$. Introduce the process $Y_t : t \mapsto e^t f(X_t)$, where

$$f(x) = \exp\left(\int_0^x (k - a(u))^{-1}\, du\right)$$

for every $x \in \mathbb{R}$. After some easy computations relying on Itô's formula (15), one can show that $Y$ verifies the SDE

$$dY_t = b(t, Y_t) \diamond dZ_t^k \qquad (18)$$

for some positive function $b$, where, as in Proposition 3, we set $Z_t^k = Z_t + kt$ for every $t \geq 0$. Fix $t = 1$ for simplicity. Notice that $f$ is a diffeomorphism, so that $Y_1 \ll \lambda$ if and only if $X_1 \ll \lambda$. Since $Z_1^k$ is obviously absolutely continuous, the problem is seemingly reduced to the case $a \equiv 0$, that is, to Proposition 3, save that $b$ also depends on $t$. However, this makes a big difference. Indeed, Fabes and Kenig [6] had constructed a continuous function $b : \mathbb{R}^2 \to [1, 2]$ such that $Y_1$ defined by (18) driven by Brownian motion is continuous singular. Notice also, still in the Brownian case, that Martini [15] had proved that $Y_1$ does not weight points as soon as $b \geq c > 0$.

Note added in proof. After this paper had been accepted, we became aware of recent related results by O. M. Kulik [in Malliavin calculus for Lévy processes with arbitrary Lévy measures. *Teor. i Imovir. Mat. Stat.* **72** (2005) 67–83]. In this paper, the author proves a general absolute continuity result for jumping SDEs without conditions on the Lévy measure (Theorem 3.1), which proves our Theorem A when Z has bounded variations.

**Acknowledgments.** We are indebted to Michel Lifshits for giving us some decisive insight in proving Theorem A, and to the referee for very careful reading. Part of this work was done during a visit to the Friedrich-Schiller-Universität in Jena by the second-named author, who would like to thank Werner Linde for his kind hospitality.

Université Henri Poincaré
Institut de Mathématiques Élie Cartan
B.P. 239
54506 Vandœuvre-lès-Nancy Cédex
France
E-mail: Ivan.Nourdin@iecn.u-nancy.fr

Université d'Évry-Val d'Essonne
Équipe d'Analyse et Probabilités
Boulevard François Mitterand
91025 Évry Cédex
France
E-mail: Thomas.Simon@maths.univ-evry.fr